\documentclass[conference]{IEEEtran}
\IEEEoverridecommandlockouts
\usepackage{cite}
\usepackage{xcolor}   
\usepackage{amsmath,amssymb,amsfonts}
\usepackage{algorithmic}
\usepackage{graphicx}
\usepackage{textcomp}
\usepackage{xcolor}
\usepackage{microtype}
\usepackage{graphicx}
\usepackage{subfigure}
\usepackage{booktabs} 
\usepackage{amsfonts}
\usepackage{amsmath}
\usepackage{url}
\usepackage{bbm}
\usepackage{algorithm}
\usepackage{amssymb,latexsym,amsthm} 
\usepackage{graphicx}
\usepackage{amsthm} 
\theoremstyle{definition}   
\newtheorem{theorem}{Theorem}   
\newtheorem{lemma}[theorem]{Lemma}
\newcommand{\expect}{\mathbb{E}}

\newcommand{\no}[1]{#1^c}

\newcommand{\eqdef}{\stackrel{\rm def}{=}}
\newcommand{\flow}{f_{\rm low}}

\DeclareMathOperator*\globmin{globmin}
\DeclareMathOperator*\globmax{globmax}

\def\BibTeX{{\rm B\kern-.05em{\sc i\kern-.025em b}\kern-.08em
    T\kern-.1667em\lower.7ex\hbox{E}\kern-.125emX}}
\begin{document}

\title{Quadratic and Cubic Regularisation Methods with Inexact function and Random Derivatives for Finite-Sum Minimisation}

\author{\IEEEauthorblockN{Stefania Bellavia}
\IEEEauthorblockA{\textit{Department of Industrial Engineering} \\
\textit{Universit\`{a} degli Studi di Firenze}\\
Firenze, Italy \\
stefania.bellavia@unifi.it}
\and
\IEEEauthorblockN{Gianmarco Gurioli}
\IEEEauthorblockA{\textit{Scuola Internazionale Superiore di}\\
\textit{Studi Avanzati (SISSA)}\\
Trieste, Italy \\
ggurioli@sissa.it}
\and
\IEEEauthorblockN{Benedetta Morini}
\IEEEauthorblockA{\textit{Department of Industrial Engineering} \\
\textit{Universit\`{a} degli Studi di Firenze}\\
Firenze, Italy \\
benedetta.morini@unifi.it}

\and
\IEEEauthorblockN{Philippe L. Toint}
\IEEEauthorblockA{\textit{Namur Center for Complex Systems (naXys)} \\
\textit{University of Namur}\\
Namur, Belgium \\
philippe.toint@unamur.be}
}

\maketitle

\begin{abstract}
This paper focuses on regularisation methods using models up to the third order to search for up to second-order critical points of a finite-sum minimisation problem.
The variant presented belongs to the framework of \cite{BellGuriMoriToin20}:
it employs random models with accuracy guaranteed with a
sufficiently large prefixed probability and deterministic inexact function
evaluations within a prescribed level of accuracy. Without assuming unbiased
estimators, the expected number of iterations is \boldmath$\mathcal{O}\bigl(\epsilon_1^{-2}\bigr)$\unboldmath or  \boldmath$\mathcal{O}\bigl(\epsilon_1^{-{3/2}}\bigr)$\unboldmath when searching for a first-order critical point using a second or third order model, respectively, and of
\boldmath$\mathcal{O}\bigl(\max[\epsilon_1^{-{3/2}},\epsilon_2^{-3}]\bigr)$\unboldmath when
seeking for second-order critical points with a third order model, in which
\boldmath$\epsilon_j$, $j\in\{1,2\}$, is the $j$th-order tolerance. \unboldmath These results
match the worst-case optimal complexity for the deterministic counterpart of
the method. Preliminary numerical tests for first-order optimality 
in the context of nonconvex binary classification in imaging, with and without Artifical Neural Networks (ANNs), are presented and discussed.
\end{abstract}


\section{Introduction}
\label{Intro} We consider adaptive regularisation methods to compute an approximate $q$-th order, $q\in\{1,2\}$, local minimum of the finite-sum minimisation problem \begin{equation}
\label{problem}
\min_{x \in \mathbb{R}^n}  f(x)\eqdef \min_{x \in \mathbb{R}^n}\sum_{i=1}^N  f_i(x), 
\end{equation}
where $f_i:\mathbb{R}^n\rightarrow\mathbb{R}$, $i\in\{1,...,N\}$. 

Problem (\ref{problem}) has recently received a large attention since it includes  a variety of applications, such as 
least-squares approximation and  
Machine Learning,  and covers  optimization problems in imaging. In fact, 
Convolutional Neural Networks (CNNs) have been  successfully employed for  classification in imaging and  
for recovering an image from a set of noisy measurements, and compares favourably with well accessed iterative reconstruction methods combined with suitable regularization. Investigating the link between traditional 
approaches  and CNNs is an  active area of research as well as the efficient solution of 
nonconvex optimization problems which arise  from training the networks on large databases of images and can be cast in the form 
(\ref{problem}) with $N$ being a large positive integer (see e.g., \cite{springer_survey, BCCPP,JMFU, MJU, shanmugamani18}).

The wide range of methods used in literature to solve 
\eqref{problem} can be classified as first-order methods,
requiring only the gradient of the objective $f$ and
second-order procedures, where the Hessian is also needed. Although
first-order schemes are generally characterised by a simple and low-cost iteration, their performance can be seriously
hindered by ill-conditioning and their success is highly dependent on the
fine-tuning of hyper-parameters. In addition, the objective function in \eqref{problem}
can be nonconvex, with a variety of local minima and/or saddle points. For
this reason, second-order strategies using curvature information have been recently used as an
instrument for escaping from saddle points more easily (see, e.g., \cite{review,CGToint,Roosta_2p,BeraBollaNoce20}). 
Clearly, the per-iteration cost is expected to be
higher than for first-order methods, since second-order derivatives
information is needed. By contrast, second-order methods have been shown to be
significantly more resilient
to ill-conditioned and badly-scaled problems, less sensitive to the choice of
hyper-parameters and tuning \cite{BeraBollaNoce20,Roosta_2p}. 

In this paper we propose methods up to order two, 
building on the Adaptive Regularization (AR)  approach.
 We stress that methods  employing cubic models show optimal complexity: 
  a first- or second-order stationary
point can be achieved in at most  $\mathcal{O}\bigl(\epsilon_1^{-{3/2}}\bigr)$
or $\mathcal{O}\bigl(\max[\epsilon_1^{-{3/2}},\epsilon_2^{-3}]\bigr)$
iterations  respectively, with $\epsilon_1$, $\epsilon_2$ being the order-dependent
requested accuracy thresholds (see, e.g., \cite{BellGuriMoriToin19,Toint1,CGToint,CGtextrmA,ARC2}). Allowing
inexactness in the function and/or derivative evaluations of AR methods,
and still preserving convergence properties and optimal complexity, has been a
challenge in recent years \cite{IMA,Cin,Roosta,CartSche17,kl,Roosta_2p,Roosta_inexact,zhou_xu_gu}. A
first approach is to impose suitable accuracy requirements that $f$ and the
derivatives have to deterministically fulfill at each iteration. But in
machine learning applications, function and derivative evaluations are generally
approximated by using uniformly randomly selected subsets of terms in the sums, and the accuracy levels can be
satisfied only within a certain probability
\cite{BellGuriMoriToin19,IMA,Cin,Roosta,Roosta_inexact,BeraBollaNoce20,CartSche17}. This suggests a
stochastic analysis of the expected
worst-case number of iterations needed to find a first- or second-order
critical point \cite{CartSche17,zhou_xu_gu}. We pursue this approach and
our contributions are as follows. We elaborate on \cite{BellGuriMoriToin20}
where adaptive regularisation methods with random models  have been proposed  
for computing strong approximate minimisers of any order for constrained smooth
optimization. We focus here on problems of the form (\ref{problem}) and on  methods from this
    class using models based on adding a regularisation term to a truncated Taylor serie of order up to two.
While gradient and Hessian are subject to random noise, function values are required to be approximated with a
deterministic level of accuracy. Our approach is  particularly suited
for applications where evaluating derivatives is more expensive than
performing function evaluations. This is for instance the case of deep neural networks
training (see, e.g., \cite{DLBook}). We discuss a matrix-free implementation in the case where gradient and Hessian are  approximated via  sampling and with adaptive accuracy requirements. The outlined procedure retains the optimal worst-case complexity and our analysis complements that in \cite{CartSche17}, as we cover the approximation of second-order minimisers.

\subsubsection{Notations} We use $\|\cdot\|$ to indicate the $2$-norm (matrices and vectors). $\nabla_x^0 f(x)$ is a synonym for $f(x)$, given $f:\mathbb{R}^n\rightarrow \mathbb{R}$. $\expect[X]$ denotes the expected value of a random variable $X$. In addition, given a random event $E$, $Pr(E)$
denotes the probability of $E$. All inexact quantities are denoted by an
overbar.

\section{Employing inexact evaluations}

\subsection{Preliminaries.} Given $1\le q\le p\le 2$, we make the following assumptions on \eqref{problem}.

\begin{description}
\item[AS.1] There exists a constant
  $\flow$ such that $f(x)\geq \flow$ for all $x \in \mathbb{R}^n$.
  \vspace{0.2cm}
\item[AS.2] $f\in \mathcal{C}^p(\mathcal{B
})$ with $\mathcal{B}$ a convex neighbourhood of $\mathbb{R}^n$. Moreover, there exists a nonnegative constant $L_p$, such that, for all $x,y\in\mathcal{B}$:\\ $ \|\nabla_x^p f(x) - \nabla_x^p f(y)\| \leq L_p \|x-y\|$.
\end{description}

\noindent
Consequently, the $p$-order Taylor expansions of $f$ centered at $x$ with
increment $s$ are   well-defined and are  given by:
$$T_{f,1}(x,s)\eqdef f(x)-\Delta T_{f,1}(x,s),$$ 
with $\Delta T_{f,1}(x,s)\eqdef   - \nabla_x f(x)^\top s$, and 
$$
T_{f,2}(x,s)\eqdef f(x)-\Delta T_{f,2}(x,s),$$
where
\begin{equation}
\label{DT-def}
\Delta T_{f,2}(x,s)\eqdef   - \nabla_x f(x)^\top s-\frac12 s^\top \nabla_x^2 f(x)s.
\end{equation}

\subsection{Optimality conditions}
To obtain complexity results  for approximating first- and second-order minimisers, we use a compact formulation to characterise the optimality conditions. As in \cite{CartGoulToin20a}, given the order of accuracy $q\in\{1,2\}$, the tolerance vector $\epsilon=\epsilon_1$, if $q=1$, or $\epsilon=(\epsilon_1,\epsilon_2)$, if $q=2$,  we say that $x\in\mathbb{R}^n$ is a $q$-th order $\epsilon$-approximate minimiser for \eqref{problem} if
\begin{equation}
\label{strong}
\phi_{f,j} (x) \leq  \frac{\epsilon_j}{j}
\textrm{ for } j=1,...,q,
\end{equation}
where 
\begin{equation}
\label{phi-def}
\phi_{f,j}(x)\eqdef f(x) - \globmin_{\|d\| \leq 1} T_{f,j}(x,d) = \globmax_{\|d\| \leq 1} \Delta T_{f,j}(x,d).
\end{equation}

The optimaliy measure  $\phi_{f,j}$ is a nonnegative (continuous) function that can be used as a measure of closeness to $q$-th order stationary points \cite{CartGoulToin20a}.  For $\epsilon_1=\epsilon_2=0$, it reduces to the known first- and second-order optimality
conditions, respectively. Indeed, assuming that $q=1$  we have from \eqref{strong}--\eqref{phi-def} with $j=1$ that $\phi_{f,1} (x)=\|\nabla_x f(x)\|=0$. If $q=2$,  \eqref{strong}--\eqref{phi-def} further imply that $\phi_{f,2}(x)= \max_{\|d\| \leq 1} \left(-\frac12 d^\top \nabla_x^2 f(x)d\right)=0$,
which is the same as requiring the semi-positive definiteness of $\nabla_x^2
f(x)$. 

\subsection{The Inexact Adaptive Regularisation ($IAR_{qp}$) algorithm}
We now define our Inexact Adaptive Regularisation ($IAR_{qp}$) scheme, whose purpose is to find a $q$-th $\epsilon$-approximate minimiser of \eqref{problem} (see \eqref{strong}) using a $p+1$ order model, for $1\le q\le p\le 2$.

The scheme, sketched in Algorithm \ref{algo} on page \pageref{algo}, is defined in analogy with the  Adaptive Regularization Approach   (see, \cite{CartGoulToin20b}), but now uses the inexact values $\overline{f}(x_k)$, $\overline{f}(x_k+s)$, $\overline{\nabla_x^j f}(x_k)$ instead of $f(x_k)$, $f(x_k+s)$, $\nabla_x^j f(x_k)$, $j=1,...,q$, respectively. At iteration $k$, the model
\begin{eqnarray}
m_k(s)&= \left\{ \begin{array}{ll}
  -\overline{\Delta T}_{f,1}(x_k,s)+\frac{\sigma_k}{2}\|s\|^2,     &\textrm{if} ~~p=1, \label{model1}\\
  -\overline{\Delta T}_{f,2}(x_k,s)+\frac{\sigma_k}{6}\|s\|^3, & \textrm{if} ~~p=2.  \label{model2}
 \end{array} \right.
\label{model} 
\end{eqnarray} 
is built at Step 1 and approximately minimised, at Step 2, finding a step $s_k$ such that
\begin{equation}
\label{descent}
m_k(s_k) \leq  m_k(0)=0
\end{equation}
 and 
$$
\overline{\phi}_{m_k,j}(s_k) =  \max_{ \|d\| \leq 1} \overline{\Delta T}_{m_k,j}(s_k,d) \leq \theta \frac{ \epsilon_j}{j}.
$$
for $j=1,...,q$, and some $\theta\in (0,\frac12)$.
Taking into account that, 
\begin{eqnarray}
\overline{\Delta T}_{m_k,1}(s_k,d)&=&-\overline{\nabla_s m}(s_k)^\top d\nonumber\\
\overline{\Delta T}_{m_k,2}(s_k,d)&=&-\overline{\nabla_s m}(s_k)^\top d-\frac12 d^\top\overline{\nabla_s^2 m}(s_k)d,\nonumber
\end{eqnarray}
we have 
\begin{eqnarray}
\overline{\phi}_{m_k,1}(s_k)&=&\max_{ \|d\| \leq 1} \left(-\overline{\nabla_s m}(s_k)^\top d\right)\nonumber \\
& =& \|\nabla_s m_k(s_k)\| \le  \theta  \epsilon_1  \label{step-term}
\end{eqnarray}
for $p=1,2$ and 
$$
\overline{\phi}_{m_k,2}(s_k)=\max_{ \|d\| \leq 1} \left( -\overline{\nabla_s m}(s_k)^\top d-\frac12 d^\top\overline{\nabla_s^2 m}(s_k)d \right)\le  \theta  \epsilon_2 
$$
 for $q=p=2$.
%


\begin{algorithm}[H]
\caption{The $IAR_{qp}$ Algorithm, $1\le q\le p\le 2$.}
\label{algo}
Step 0: Initialization.
\begin{description}
\item[]
An initial point $x_0\in\mathbb{R}^n$, initial regulariser
  $\sigma_0>0$, tolerances $\epsilon_j$, $j=1,...,q$, and constants $\theta\in(0,\frac12)$, $\eta \in (0,1)$,  $\gamma>1$,
  $\alpha \in (0,1)$, $\omega_0=\min \left[\frac12 \alpha \eta,\frac{1}{\sigma_0}\right]$, $\sigma_{\min}\in (0,\sigma_0)$ are given. Set $k=0$.
\end{description}
Step 1: Model definition. 
\begin{description}
\item[] Build the approximate gradient $\overline{ \nabla_x f}(x_k)$. If $p=2$, compute also the Hessian $\overline{\nabla_x^2 f}(x_k)$. Compute the model
  $m_k(s)$ as defined in \eqref{model}.
\end{description}
Step 2: Step calculation. 
\begin{description}
\item[]
  Compute a step $s_k$  satisfying \eqref{descent}--\eqref{step-term}, for $j=1,...,q$. If $\overline{\Delta T}_{f,p}(x_k,s_k) = 0$, go to Step~4.
\end{description}
Step 3: Function approximations.
\begin{description}
\item[ ] 
  Compute $\overline{f}(x_k)$ and $\overline{f}(x_k+s_k)$ satisfying
  \eqref{f-est1}--\eqref{f-est2}.
  \end{description}
Step 4: Test of acceptance.
  \begin{description}
\item[ ] Set
\[
  \rho_k = \left\{\begin{array}{ll}
  \frac{\overline{f}(x_k) - \overline{f}(x_k+s_k)}{\overline{\Delta T}_{f,p}(x_k,s_k)} &
  \textrm{ if ~} \overline{\Delta T}_{f,p}(x_k,s_k) > 0, \\
  -\infty & \textrm{ otherwise.}
  \end{array}\right.
\]

  If $\rho_k \geq \eta$ (\textit{successful iteration}), then set
  $x_{k+1} = x_k + s_k$; otherwise (\textit{unsuccessful iteration}) set $x_{k+1} = x_k$.
\end{description}
Step 5: Regularisation parameter update.
\begin{description}  
\item[ ]
  Set
 \begin{equation}
 \label{sigupdate}
 \sigma_{k+1} = \left\{ \begin{array}{ll}
   \max\left[\sigma_{\min},\frac{1}{\gamma} \sigma_k\right],     &\textrm{if} \rho_k \geq \eta,\\
   \gamma \sigma_k, & \textrm{if} \rho_k < \eta.
 \end{array} \right.
 \end{equation} 
 \end{description}
Step 6: Relative accuracy update.
 \begin{description}
\item [ ]
 Set
\begin{equation}
\label{omega-def}
 \omega_{k+1} = \min \left[\frac12 \alpha \eta,\frac{1}{\sigma_{k+1}}\right].
\end{equation}
 Increment $k$ by one and go to Step~1.
\end{description}
\end{algorithm}

The existence of such a step can be proved as in Lemma~$4.4$ of
\cite{CartGoulToin20a} (see also \cite{BellGuriMoriToin20}). We note  that the model definition in \eqref{model} does not
depend on the approximate value of $f$ at $x_k$.  The ratio $\rho_k$, depending on inexact function values and model,  is then computed at Step $4$ and affects the acceptance of the
trial point $x_k+s_k$. Its magnitude also influences the regularisation
parameter $\sigma_k$ update at Step $5$. While the gradient
$\overline{\nabla_x f}(x_k)$ and the Hessian $\overline{\nabla_x^2 f}(x_k)$
approximations can be seen as random estimates, the values
$\overline{f}(x_k)$, $\overline{f}(x_k+s_k)$ are required to be
deterministically computed to satisfy
\begin{eqnarray}
\left|\overline{f}(x_k)-f(x_k)\right|&\le& \omega_k \overline{\Delta T}_{f,p}(x_k,s_k),~~~~~~ \label{f-est1}\\
\left|\overline{f}(x_k+s_k)-f(x_k+s_k)\right|&\le& \omega_k \overline{\Delta T}_{f,p}(x_k,s_k),~~~~~~ \label{f-est2}
\end{eqnarray}
in which $\omega_k$ is iteratively defined at Step~$6$. As for the implementation of the  algorithm, 
we note that 
$\overline{\phi}_{m_k,1}(s_k)=\|\nabla_s m_k(s_k)\|$,
while $\overline{\phi}_{m_k,2}(s_k)$ can be computed via a
standard trust-region method at a cost which is comparable to that of
computing the Hessian left-most eigenvalue.
The approximate minimisation of the cubic model \eqref{model} at each
iteration can be seen, for $p=2$, as an issue in the AR framework. However, an approximate minimiser can be computed via matrix-free approaches 
accessing the Hessian only
through matrix-vector products.
A number of procedures have been proposed, 
ranging from
Lanczos-type iterations where the minimisation is done via nested, lower
dimensional, Krylov subspaces \cite{ARC1}, up to minimisation via 
gradient descent (see, e.g.,
\cite{Agar16,CarmDuch16,CarmDuchHindSidf18}) or the  Barzilai-Borwein gradient method \cite{BLMS15}. 
 Hessian-vector products  can be approximated by the finite difference approximation,
with only two gradient evaluations
\cite{BLMS15,CarmDuchHindSidf18}. All these matrix-free implementations remain
relevant if $\nabla_x^2 f(x_k)$ is  defined via subsampling, proceeding as in Section \ref{Section_PNT}.   Interestingly,
back-propagation-like methods in deep learning also allow computations of
Hessian-vector products at a similar cost \cite{Pear94,Schr02}.

\subsection{Probabilistic assumptions on $IAR_{qp}$}
In what follows, all
random quantities are denoted by capital letters, while the use of small
letters denotes their realisations. We refer to the random model
$M_k$ at iteration $k$, while $m_k$ = $M_k(\zeta_k)$ is its realisation, with
$\zeta_k$ being a random sample taken from a context-dependent probability space. As
a consequence, the iterates $X_k$, as well as the regularisers $\Sigma_k$, the
steps $S_k$ and $\Omega_k$, are the random variables such that
$x_k=X_k(\zeta_k)$, $\sigma_k=\Sigma_k(\zeta_k)$, $s_k=S_k(\zeta_k)$ and
$\omega_k=\Omega_k(\zeta_k)$. For the sake of brevity, we will
omit $\zeta_k$ in what follows. Due to the randomness of the model construction at Step $1$, the $IAR_{qp}$ algorithm induces a random process formalised by $\{X_k, S_k, M_k, \Sigma_k,\Omega_k\}$ ($x_0$ and $\sigma_0$ are deterministic quantities). For $k\ge 0$, we formalise the conditioning on the past by using
$\mathcal{A}_{k-1}^{M}$, the $\hat{\sigma}$-algebra induced by the random
variables $M_0$, $M_1$,..., $M_{k-1}$, with
$\mathcal{A}_{-1}^{M}=\hat{\sigma}(x_0)$. We also denote by $d_{k,j}$ and
$\overline{d}_{k,j}$ the arguments in the maximum in the definitions of
$\phi_{m_k,j}(s_k)$ and
$\overline{\phi}_{m_k,j}(s_k)$, respectively. We say that iteration $k$ is \textit{accurate} if
 \begin{equation}
\label{dis1}
  \|\overline{\nabla_x^\ell f}(X_k) -\nabla_x^\ell f(X_k)\|
  \leq \Omega_k \frac{\overline{\Delta T}_{k,\min}}{6\tau_k^\ell},
  \textrm{ for }\ell\in\{1,...,p\},
\end{equation}
where
\begin{eqnarray}
\tau_k &\eqdef& \max\left[ \|S_k\|, \max_{j=1,...,q}[\|D_{k,j}\|, \|\overline{D}_{k,j}\|]\right]\nonumber\\
 \overline{\Delta T}_{k,\min} &\eqdef& \min \Bigl[\overline{\Delta T}_{f,p}(X_k,S_k),\nonumber\\
    && ~~~~~~\min_{j=1,...,q}\Bigl[ \overline{\Delta T}_{m_k,j}(S_k,D_{k,j}),\nonumber\\
    && ~~~~~~~~~~~~~~~~~  \overline{\Delta T}_{m_k,j}(S_k,\overline{D}_{k,j})\Bigr]\Bigr].\nonumber
 \end{eqnarray}
We emphasize that the above accuracy requirements are adaptive.
At variance with the trust-region
methods of \cite{BlanCartMeniSche19,CartSche17,ChenMeniSche18}, the above conditions do not require the model to be fully linear or quadratic in a ball centered at $x_k$ of radius at least $\|s_k\|$.  
As standard in related papers \cite{BlanCartMeniSche19,CartSche17,ChenMeniSche18,PaquSche20}, we assume a lower bound on the probability of the model to be accurate at the $k$-th iteration.
\begin{description}
\item[AS.3] For all $k\ge 0$, conditioned to $\mathcal{A}_{k-1}^{M}$, we assume that \eqref{dis1} is satisfied with probability at least $p_* \in (\frac12,1]$ independent of $k$.
 \end{description}
We stress that the inequalities in
\eqref{dis1} can be satisfied via uniform subsampling with
probability of success at least $(1-t_{\ell})$, using the operator Bernstein inequality (see, Subsection \ref{Sub_II} and Section $7.2$ in \cite{BellGuriMoriToin19}), and this provides AS.3 with $p_*=\prod_{\ell=1}^p (1-t_{\ell})$.

\section{Worst-case complexity analysis} 

\subsection{Stopping time}

Before starting with the worst-case evaluation complexity of the $IAR_{qp}$
algorithm, the following clarifications are important. For each
$k\ge 1$ we assume that the computation of $s_{k-1}$ and thus of the trial
point $x_{k-1}+s_{k-1}$ ($k\ge 1$) are deterministic, once the inexact model
$m_{k-1}(s)$ is known, and that \eqref{f-est1}-\eqref{f-est2} at Step $3$ of
the algorithm are enforced deterministically; therefore, $\rho_{k-1}$ and the
fact that iteration $k-1$ is successful are deterministic outcomes of the
realisation of the (random) inexact model. Hence, 
\[
N_\epsilon
=\inf \Bigl\{ k\geq 0~\mid~ \phi_{f,j} (X_k) \leq \frac{\epsilon_j}
{j},~ j =1,...,q\Bigr\}
\]
can be seen as a family of hitting times depending of $\epsilon$ and corresponding to the number of iterations required until \eqref{strong} is met for the first time. The assumptions made on top of this subsection
imply that the variables $X_{k-1}+S_{k-1}$ and the event $\{X_k=X_{k-1}+S_{k-1}\}$, occurring when iteration $k-1$ is successful, are measurable with respect to $\mathcal{A}_{k-1}^M$. Our aim is then to derive an upper bound on the expected number of steps $\expect[N_\epsilon]$ needed by the $IAR_{qp}$ algorithm, in the worst-case, to reach an $\epsilon$-approximate $q$-th-order-necessary minimiser, as in \eqref{strong}.

\subsection{Deriving expected upper bounds on $N_{\epsilon}$}A crucial property for our
complexity analysis is to show that, when the model \eqref{model} is accurate,
iteration $k$ is successful but $\|\nabla_x f(x_{k+1})\|>\epsilon_1$, then
$\|s_k\|^{p+1}$ is bounded below by $\psi(\sigma_k)\epsilon_1^{\frac{p+1}{p-q+1}}$ in which
$\psi(\sigma)$ is a decreasing function of $\sigma$. This is
central in \cite{CartSche17} for proving the
worst-case complexity bound for first-order optimality. By virtue of the compact
formulation \eqref{phi-def} of the optimality measure, for $p=2$, the lower bound on $\|s_k\|^3$
also holds for second-order optimality and a suitable power of $\epsilon_2$.

\begin{lemma}
Suppose that AS.2 holds and consider any realisation of the
  algorithm. Suppose that \eqref{dis1} also occurs, that iteration $k$ is
  successful and that, for some $j=1,...,q$, \eqref{strong} fails
  for $x_{k+1}$. Then 
  $\|s_k\|^{p+1} \geq \psi(\sigma_k) \epsilon_j^{\frac{p+1}{p-q+1}}$,
  with
  \begin{equation}
    \label{psi-def1}
 \psi(\sigma) =  \min\left[ 1,\left(\frac{(1-2\theta)(3-q)}
     {q(L_{f,2}+\sigma)}\right)^{\frac{p+1}{p-q+1}}\right].
  \end{equation}
\end{lemma}
\noindent
We refer to Lemma $3.4$ in \cite{BellGuriMoriToin20} for a complete proof. In order to develop the complexity analysis, building on \cite{CartSche17} and
using results from the previous subsection, we first identify a
threshold $\sigma_s$ for the regulariser, above which each iteration of the
algorithm with accurate model is successful. We also give some preliminary bounds (see \cite{BellGuriMoriToin20}).

\begin{lemma}
\label{acc-sigma-success-l}
Suppose that AS.2 holds. For any realisation of the algorithm, if iteration
$k$ is such that \eqref{dis1} occurs and 
\begin{equation}
\label{large-sigma}
\sigma_k \geq \sigma_s \eqdef \max\left[2\sigma_0,\frac{L_{f,p}+3}{1-\eta}\right],
\end{equation}
then iteration $k$ is successful.
\end{lemma}

\begin{lemma}
Let Assumptions AS.1--AS.3 hold. For all realisations of the $IAR_{qp}$, let $N_{AS}$, $N_{\Lambda}$ and $N_{\Lambda^c}$ represent the number of accurate successful iterations with $\Sigma_k \leq \sigma_s$, the number of iterations with   $\Sigma_k<\sigma_s$ and the number of iterations with  $\Sigma_k \ge \sigma_s$, respectively. Then,
\begin{eqnarray}
\expect[N_{\no{\Lambda}}] &\leq& \frac{1}{2p_*}\expect[N_\epsilon],\nonumber\\
 \expect[N_{AS}]&\le& \frac{6(f_0-\flow)}{\eta\sigma_{\min}(1-\alpha)\psi(\sigma_s)}\max_{j=1,...,q}\left[\epsilon_j^{-\frac{p+1}{p-q+1}}\right]+1,\nonumber\\
 \expect[N_{\Lambda}]&\leq& \frac{1}{2p_*-1}\left[2\expect[N_{AS}]+ \log_{\gamma}\Bigl( \frac{\sigma_s}{\sigma_0}  \Bigr)\right].\nonumber
\end{eqnarray}
\end{lemma}
\noindent
The proofs of the first and the third bound of the previous lemma follow the
reasoning of \cite{CartSche17}, the proof of the second bound is given in
\cite{BellGuriMoriToin20}. Finally, the fact that
$\expect[N_{\epsilon}]=\expect[N_{\Lambda}]+\expect[N_{\no{\Lambda}}]$ allows
us to state our complexity result.

\begin{theorem}
\label{ThComplexity}
Under the assumptions of Lemma $3$,
\begin{eqnarray}
\expect[N_\epsilon]
&\leq& \frac{2p_*}{(2p_*-1)^2}
\Bigl[\frac{12(f_0-\flow)}{\eta\sigma_{\min}(1-\alpha)\psi(\sigma_{s})}
  \max_{j=1,...,q}\left[\epsilon_j^{-\frac{p+1}{p-q+1}}\right]
  \nonumber\\
  &&~+\log_{\gamma}\left( \frac{\sigma_s}{\sigma_0}\right)+2 \Bigr].\label{final-bound}
\end{eqnarray}
\end{theorem}

Theorem \ref{ThComplexity} shows that for $q=p=1$ the $IAR_{1}$ algorithm reaches the hitting time in at most $\mathcal{O}(\epsilon_1^{-2})$ iterations. Thus, it matches the complexity bound of gradient-type methods for nonconvex problems employing exact gradients. Moreover, $IAR_{qp}$ needs  $\mathcal{O}(\epsilon_1^{-3/2})$ iterations when $q=1$ and $p=2$ are considered, while to approximate a second-order optimality point (case $p=q=2$) $\mathcal{O}\bigl(\max[\epsilon_1^{-{3/2}},\epsilon_2^{-3}]\bigr)$ are needed in the worst-case. Therefore, for $1\le q\le p \le 2$, Theorem \ref{ThComplexity} generalises the complexity bounds stated in \cite[Theorem 5.5]{CartGoulToin20a} to the case where evaluations of $f$ and its derivatives are inexact, under probabilistic assumptions on the accuracies of the latter. Since in  \cite[Theorems 6.1 and 6.4]{CartGoulToin20a} it is shown that the complexity bounds are sharp in order of the tolerances $\epsilon_1$ and $\epsilon_2$, for exact evaluations and Lipschitz continuous derivatives of f, this is also the case for the bounds of Theorem \ref{ThComplexity}.

{{
\section{Preliminary numerical tests}\label{Section_PNT}

We here present preliminary novel numerical results on the $IAR_{qp}$ algorithm with model up to order three ($p\in\{1,2\}$) and first-order optimality ($q=1$). Numerical tests are performed within the context of nonconvex finite-sum minimisation problems for the binary digits detection of two datasets coming from the imaging community. Specifically, we consider the MNIST dataset in \cite{Mnist}, usually considered for classifying the $10$ handwritten digits and renamed here as MNIST-B and used to discard even digits from odd digits, and the GISETTE database \cite{UCI}, to detect the highly confusable handwritten digits ``4'' and ``9''.\\
The overall binary classification learning procedure focuses on the two macro-steps reported below.
\begin{itemize}
\item \textit{The training phase}. Given a training set $\{(a_i,y_i)\}_{i=1}^N$ of $N$ features  $a_i\in\mathbb{R}^d$  and corresponding binary labels $y_i\in\{0,1\}$,
we solve the following minimisation problem:
\begin{equation}
\label{minloss}
\begin{split}
\min_{x\in\mathbb{R}^n} f(x)&= \min_{x\in\mathbb{R}^n}  \frac{1}{N}\sum_{i=1}^N{f_i(x)}\\
&\eqdef \min_{x\in\mathbb{R}^n}\frac 1 N \sum_{i=1}^N \left( y_i-net\left(a_i; x\right) \right)^2,
\end{split}
\end{equation}
in which the objective function is the so-called \textit{training loss}.
We use a feed-forward Artificial Neural Network (ANN) $net: a\in \mathbb{R}^d \rightarrow (0,1)$, defined by the parameter vector $x\in\mathbb{R}^n$, 
as the model for predicting the values of the labels, and the square loss as a measure of the error on such predictions, that has to be minimised by approximately solving \eqref{minloss} in order to come out with the parameter vector $x^*$, to be used for label predictions on the testing set of the following step. We use the notation $(d_1, ...., d_h)$ to denote a feed-forward ANN with $h$ hidden layers where the $j$th hidden layer, $j\in\{1,...,h\}$, has $d_j$ neurons. The number of neurons for the input layer and the output layer are constrained to be $d$ and $1$, respectively. The $tanh$ is considered as the activation function for each inner layers, while the sigmoid $\sigma$ constitutes the activation at the output layer. Therefore, if zero bias and no hidden layers ($h=0$) are considered, $d=n$ and $net(a; x)$ reduces to $\sigma(a^\top x)=1/(1+e^{-a^\top x})$, that constitutes our no net model.

\item \textit{The testing phase}. Once the training phase is completed, a number $N_T$ of testing data $\{\overline a_i,\overline y_i\}_{i=1}^{N_T}$ is used to validate the computed model. The values $net(\overline{a}_i; x^*)$ are used to predict the testing labels $\overline y_i$, $i\in\{1,...,N_T\}$, computing the corresponding error (\textit{testing loss}), measured by
$
\frac{1}{N_T} \sum_{i=1}^{N_T} \left(\overline y_i-net(\overline{a}_i; x^*) \right)^2.
$
\end{itemize}

\subsection{Implementation issues and results}\label{Sub_II}
The implementation of the $IAR_{qp}$ algorithm respects the following specifications. The cubic regularisation parameter is initialised as $\sigma_0=10^{-1}$, its minimum value is $\sigma_{\min}=10^{-5}$ and the initial guess vector $x_0=(0,...,0)^\top \in\mathbb{R}^n$ is considered for all runs. Moreover, the parameters $\alpha=0.5$,  $\eta=0.8$ and $\gamma=2$ are used. }
The estimators of the true objective functions and derivatives are obtained by averaging in a subsample of $\{1,\ldots,N\}$. More specifically, the approximations of the objective function values and of first and second
derivatives take the form:
\begin{equation}\label{approxf}
\begin{split}
&\overline{f}_k (x_k)= \frac{1}{|{\cal D}_{k,1}|} \sum_{i \in {\cal D}_{k,1}} f_i(x_k),\\
&\overline{f}_k (x_k+s_k)= \frac{1}{|{\cal D}_{k,2}|} \sum_{i \in {\cal D}_{k,2}} f_i(x_k+s_k),
\end{split}
\end{equation}
\begin{eqnarray}
\overline{\nabla f}(x_k)
= \frac{1}{|{\cal G}_k|} \sum_{i \in {\cal G}_k} \overline{\nabla f_i}(x_k),\label{approxg}\\
\overline{\nabla^2 f}(x_k)
= \frac{1}{|{\cal H}_k|} \sum_{i \in {\cal H}_k} \overline{\nabla^2 f_i}(x_k), \label{approxH}
\end{eqnarray}
\noindent
where ${\cal D}_{k,1}$, ${\cal D}_{k,2}$, ${\cal G}_k$, ${\cal H}_k$ are subsets of
$\{1,2,\ldots,N\}$ with sample sizes $|\mathcal{D}_{k,1}|, |\mathcal{D}_{k,2}|, |\mathcal{G}_k|,|{\cal H}_k|$.
These  sample sizes are adaptively chosen by the procedure in order to satisfy \eqref{f-est1}--\eqref{dis1} with probability at least $1-t$. We underline that the enforcement of  \eqref{f-est1}--\eqref{f-est2} on function evaluations is relaxed in our experiments and that such accuracy requirements are now supposed to hold in probability.
Given the prefixed absolute accuracies  $\nu_{\ell,k}$ for the derivative of order $\ell$ at iteration $k$ and $t$ a prescribed probability of failure, 
by using results in \cite{Trop15} and  \cite{BellGuriMoriToin19}, the approximations given in \eqref{approxf}--\eqref{approxH}
satisfy ($\ell\in\{1,2\}$):
\begin{eqnarray}
Pr\Big[\|\overline{\nabla^{\ell} f}(x_k)-\nabla^{\ell}f(x_k)\|\leq\nu_{\ell,k}\Big]
& \ge & 1-t,\label{vareps-j-prob}\\
  Pr\Big[|\overline{f}_k(x_k+s_k)-f(x_k+s_k)| \leq \nu_{0,k}\Big]& \ge & 1-t,\label{Df+-DT-prob}\\
Pr\Big[|\overline{f}_k(x_k)-f(x_k)| \leq \nu_{0,k}\Big] & \ge & 1-t,\label{Df-DT-prob}
  \end{eqnarray}
provided that 
\begin{equation}\label{Dkl}
 |\mathcal{D}_{k,j}|\geq\min \left \{N, \left\lceil
 \frac{4\kappa_{f,0}(x)}{\nu_{0,k}}
  \left(\frac{2\kappa_{f,0}(x)}{\nu_{0,k}}+\frac{1}{3}\right)
  \log\left(\frac{2}{t}\right)\right \rceil\right \},
\end{equation}
\[
x= \left\{ \begin{array}{ll}
 {}x_k, & \mbox{if}~ j=1, \\
 {}x_k+s_k,       &\mbox{if}~ j=2,
 \end{array} \right.
\]
\begin{equation}\label{Gk}
|\mathcal{G}_{k}|\geq\min \left \{N,\left\lceil
 \frac{4\kappa_{f,1}(x_k)}{\nu_{1,k}}
  \left(\frac{2\kappa_{f,1}(x_k)}{\nu_{1,k}}+\frac{1}{3}\right)
  \log\left(\frac{n+1}{t}\right)\right \rceil\right\},\
\end{equation}
\begin{equation}\label{Hk}
|\mathcal{H}_k| \geq\min\left \{ N, \left \lceil
 \frac{4\kappa_{f,2}(x_k)}{\nu_{2,k}}
  \left(\frac{2\kappa_{f,2}(x_k)}{\nu_{2,k}}+\frac{1}{3}\right)
  \log\left(\frac{2n}{t}\right)\right \rceil\right \}.
\end{equation}
The nonnegative constants $\{\kappa_{f,\ell}\}_{\ell=0}^2$ in \eqref{approxf}--\eqref{approxH} should be such that
(see, e.g., \cite{BellGuriMoriToin19}), for $x\in \mathbb{R}^n$ and all $\ell\in\{0,1,2\}$, $\max_{i \in\{1,...,N\}}\|\nabla^{\ell}f_i(x)\| \leq \kappa_{f,\ell}(x)$. Since their estimations can be challenging, we consider here a constant $\kappa\eqdef\kappa_{f,\ell} $, for all $\ell\in\{0,1,2\}$, setting its value experimentally, in order to control the growth of the sample sizes \eqref{Dkl}--\eqref{Hk} throughout the running of the algorithm. In particular, $\kappa=8\cdot 10^{-4}$ and $\kappa=3\cdot 10^{-2}$ are considered by $IAR_1$ on the GISETTE and MNIST-B dataset, respectively.  Moreover, in \eqref{Dkl}-\eqref{Hk}, we used $t=0.2$.

We stress that \eqref{f-est1}-\eqref{dis1} ensure that $\nu_{0,k}=\omega_k\overline{\Delta T}_{f,p}(x_k,s_k)$, $\nu_{\ell,k}=\Omega_k \frac{\overline{\Delta T}_{k,\min}}{6\tau_k^\ell},
  \textrm{ for }{{\ell\in\{1,...,p\}}}.$   
  The computation of both  $\overline{\Delta T}_{f,p}(x_k,s_k)$ and $\overline{\Delta T}_{k,\min}$ requires the knowledge of the step $s_k$, that is available at Step 3 when function estimators are built, but it is not yet available at Step 1 when approximations of the gradient (and the Hessian) have to be computed.
Thus, practical implementation of Step 1  calls for an inner loop  where the accuracy requirements is progressively reduced in order to meet  \eqref{dis1}.  This process terminates as for sufficiently small accuracy requirements the right-hand side in \eqref{approxg}-\eqref{approxH} reaches the full sample $N$. To make clear this point 
 a detailed description of Step 1 and Step 2  of the $IAR_{qp}$ algorithm with $q=p=1$, renamed for brevity as $IAR_{1}$,  is given below.
 
 \begin{algorithm}[H]
\caption{Detailed Steps 1--2 of the $IAR_1$ algorithm.}
Step 1: Model definition.
\begin{description}
\item[ ]

  Initialise $\varepsilon_{1,0} = \kappa_\varepsilon$ and set $i=0$.
  Do
  \begin{enumerate}
    \item compute $\overline{\nabla f}(x_k)$ with
      $\|\overline{\nabla f}(x_k)-\nabla f(x_k)\| \leq
      \varepsilon_{1,i}$ and increment $i$ by one.
    \item if $\varepsilon_{1,i}\leq \omega \|\overline{\nabla f}(x_k)\|$,
      go to Step 2;
    \item set $\varepsilon_{1,i+1}= \gamma_\varepsilon\varepsilon_{1,i}$ and
      return to item 1 in this enumeration.
    \end{enumerate}
  \end{description}
  Step 2: Step calculation.
 \begin{description}   
\item[ ]
  Set
  \[
  s_k = -\overline{\nabla f}(x_k)/\sigma_k,\quad
  \overline{\Delta T}_{f,1}(x_k,s_k)= \frac{\|\overline{\nabla f}(x_k)\|^2}{\sigma_k}.
  \]
\end{description}
\end{algorithm}

Regarding algorithm    $IAR_{qp}$ with $q=1$ and $p=2$, hereafter called $IAR_{2}$, a similar loop as in Step 1 of Algorithm2  is needed that involves the step computation.  The approximated gradient and Hessian  are  computed via \eqref{approxg}
 and  \eqref{approxH}, 
using  predicted accuracy requirements $\varepsilon_{1,0}$ and $\varepsilon_{2,0}$, respectively. Then, the step $s_k$ is computed and 
if  \eqref{dis1} is not satisfied 
 then the accuracy requirements 
are progressively decreased and the step is recomputed  until  \eqref{dis1} is satisfied for the first time.
 The values $\kappa_{\epsilon}=\gamma_{\epsilon}=0.5$ are considered to initialize and decrease the the accuracy requirements in both $IAR_1$ and $IAR_2$ Algorithms.  

As one may easily see,   the step $s_k$ taken in $IAR_1$ Algorithm  (see the form of the model \eqref{model1} and Step 2 in Algorithm 2) is simply  the global minimiser the model. 
Therefore, for any successful iteration $k$:
    \begin{equation}\label{AR1DAOx}
    x_{k+1}=x_k+s_k=x_k-\frac{1}{\sigma_k}\overline{\nabla f}(x_k),
    \end{equation}
    that can be viewed as the general iteration of a gradient descent method under inexact gradient evaluations with an adaptive choice of the learning rate (it corresponds to the opposite of the coefficient of $\overline{\nabla f}(x_k)$ in \eqref{AR1DAOx}), taken as the inverse of the quadratic regulariser $\sigma_k>0$.
    
On the other hand, in the $IAR_2$ procedure,
an approximated minimizer $s_k$ of the  third model $m_k$ \eqref{model2} has to be computed. 
Such approximate minimizer has to satisfy 
$m_k(s_k) \leq  m_k(0)$ and  \eqref{step-term}, i.e.
$\|\nabla_s m_k(s_k)\| \leq \theta { \epsilon_1}$. In our implementation, we used  $\theta=0.5$, $\epsilon_1=10^{-3}$ and the  computation of $s_k$
has been carried out using the Barzilai-Borwein 
method \cite{raydan1997barzilai} with a non-monotone 
line-search following \cite[Algorithm 1]{Serafino2018} and using parameter 
values as suggested therein.
The Hessian-vector product required at each iteration of the Barzilai-Borwein procedure is approximated   via finite-differences in the model gradient.

To assess the computational cost of the $IAR_{1}$ and $IAR_2$ Algorithms, we follow the approach in \cite{Roosta_2p} for what concerns the cost measure definition.
Specifically, at the generic iteration $k$, we count  the $N$  forward propagations needed to evaluate the objective in \eqref{minloss} at $x_k$ has a unit Cost Measure (CM), while the evaluation of the  approximated gradient at the same point requires an extra cost of $\left( 2|\mathcal{G}_k\smallsetminus (\mathcal{G}_k\cap\mathcal{D}_{k,1})|+ |\mathcal{G}_k\cap \mathcal{D}_{k,1}|\right)/N$ CM. Moreover,   at each iteration of the Barzilai-Borwein procedure of the $IAR_2$ Algorithm, the approximation of the Hessian-vector product by finite differences calls for the computation of
$\overline{\nabla f} (x_k+hv)$, ($h\in\mathbb{R}^+$),  leading to  additional  $|\mathcal{H}_k|$ forward and backward propagations, at the price of the weighted cost $2|\mathcal{H}_k|/N$ CM, and a potential extra-cost of $|\mathcal{H}_k\smallsetminus(\mathcal{G}_k\cap \mathcal{H}_k)|/N$  CM to 
compute $\nabla f_i(x_k)$, for $i\in\mathcal{H}_k\smallsetminus(\mathcal{G}_k\cap \mathcal{H}_k)$. A budget of $80$ and $100$ CM have been considered as the stopping rule for our runs for the MNIST-B and the GISETTE dataset, respectively.\\
We now solve the binary classification task described at the beginning of Section \ref{Section_PNT} on the MNIST-B and GISETTE datasets, described in Table \ref{TableMLData}. We test the $IAR_1$ algorithm on both datasets without net, while different ANN structures are considered for the MNIST-B via $IAR_1$. A comparison between $IAR_{2}$ and $IAR_1$ is performed on the GISETTE dataset. The accuracy in terms of the binary classification rate on the testing set, for each of the performed methods and ANN structures, is reported in Table \ref{TableANN}, averaging over 20 runs.

\begin{small}
\begin{table}[h] 
\begin{center}
\begin{tabular}{l|ccccc}
\toprule
Dataset & Training~$N$ & $d$  & Testing $N_T$ \\  \hline
\midrule

GISETTE &   $4800$ & $5000$ & $1200$\\

MNIST-B &   $60000$ & $784$ & $10000$\\
\bottomrule
\end{tabular}
\caption{MNIST-B and GISETTE datasets. Size of the training set ($N$), problem dimension ($d$) and size of the testing set ($N_T$).}\label{TableMLData}
\end{center}
\end{table}
\end{small}

\begin{small}
\begin{table}[h] 
\begin{center}
\begin{tabular}{l|c|c|c|cc}
\toprule
Dataset & $IAR_1$ & $IAR_1$ & $IAR_1$  & $IAR_2$ \\
&\mbox{(no net)} & shallow net: $(15)$ & 2-hidden layers & (no net) \\
&                         &                                  &     net: $(15,2)$   & \\  \hline
\midrule
GISETTE & $87.75\%$ & -- &   --  & $94.67\%$ \\
MNIST-B & $87.37\%$ & $88.42\%$ &   $89.23\%$  & -- \\
\bottomrule
\end{tabular}
\caption{MNIST-B and GISETTE datasets. List of the considered methods and architectures with the binary classification rate on the testing set; mean values over $20$ runs.}\label{TableANN}
\end{center}
\end{table}
\end{small}
\noindent
We immediately note that, for the series of tests performed using $IAR_1$, adding hidden layers with a small amount of neurons seems to be effective, since the accuracy in the binary classification on the testing set increases. In order to avoid overfitting, no net is considered on the GISETTE dataset, because of the fact that the number of parameters $n=d$ is greater than the number of  training samples $N_T$ (see Table \ref{TableMLData}) already in the case of no hidden layers.  We further notice that, the second order procedure provides, for the GISETTE dataset, a higher classification accuracy. Therefore, on this test, adding second order information is beneficial.  To get more insight, in Figures \ref{PerfAR1DA}--\ref{PerfANNAR1DA} and Figure \ref{PerfAR2DA} we plot  the training and testing loss against CM, while Figures \ref{SSizeAR1DA}--\ref{SSizeANNAR1DA} and Figure \ref{SSizeAR2DA} are reserved to the computed sample sizes, performed by  the tested procedures. 
Concerning the MNIST-B dataset, Figures \ref{PerfAR1DA}--\ref{PerfANNAR1DA} 
show that the decrease of the training and testing loss against CM performed by $IAR_1$ goes down to more or less the same values at termination, meaning that the training phase has been particularly effective and thus the method is able of generalising the ability of predicting labels learnt on the training phase when the testing data are considered. 
 \begin{figure}[h]
\centering
\includegraphics[width=%
0.24\textwidth]{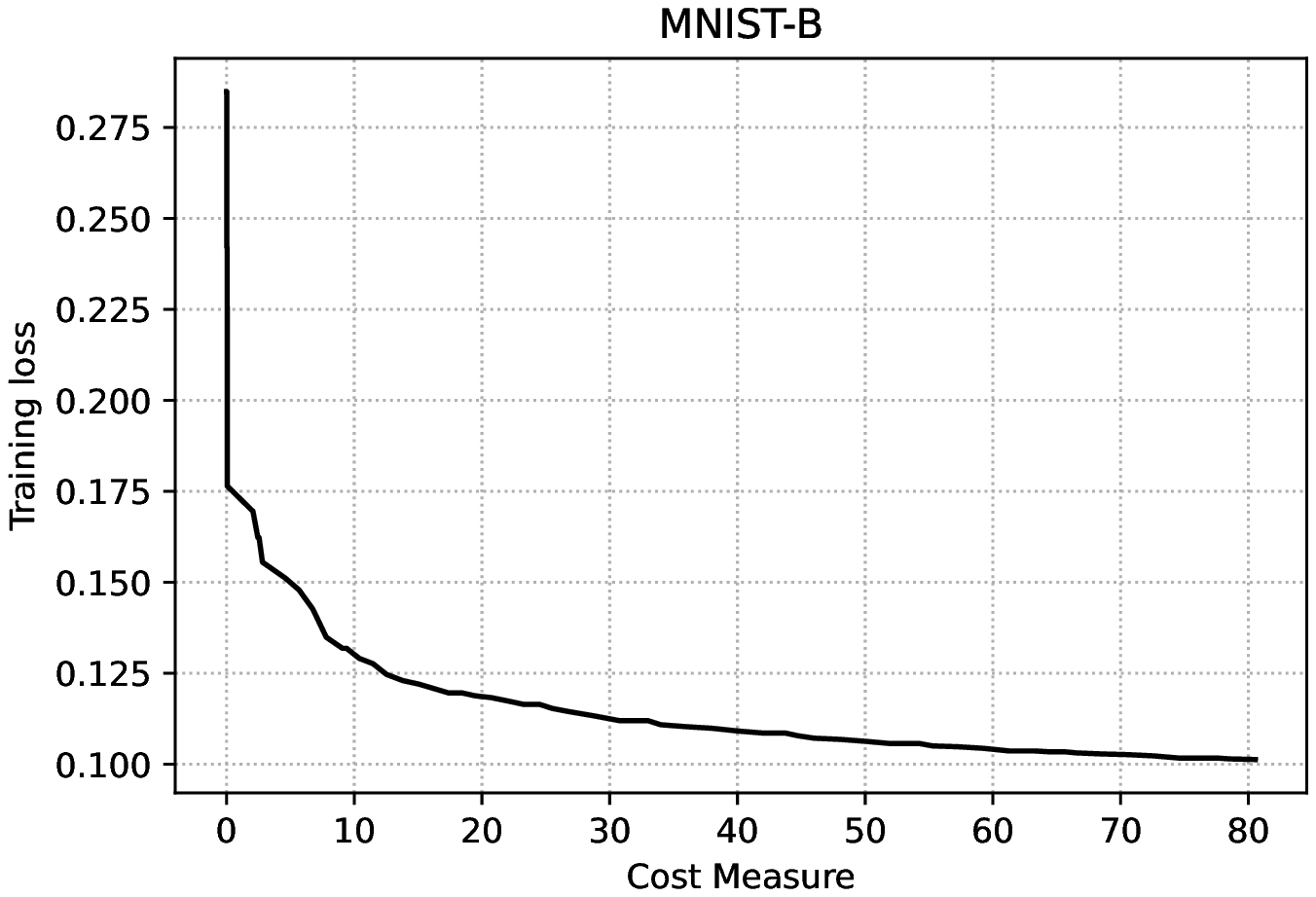}
\includegraphics[width=%
0.24\textwidth]{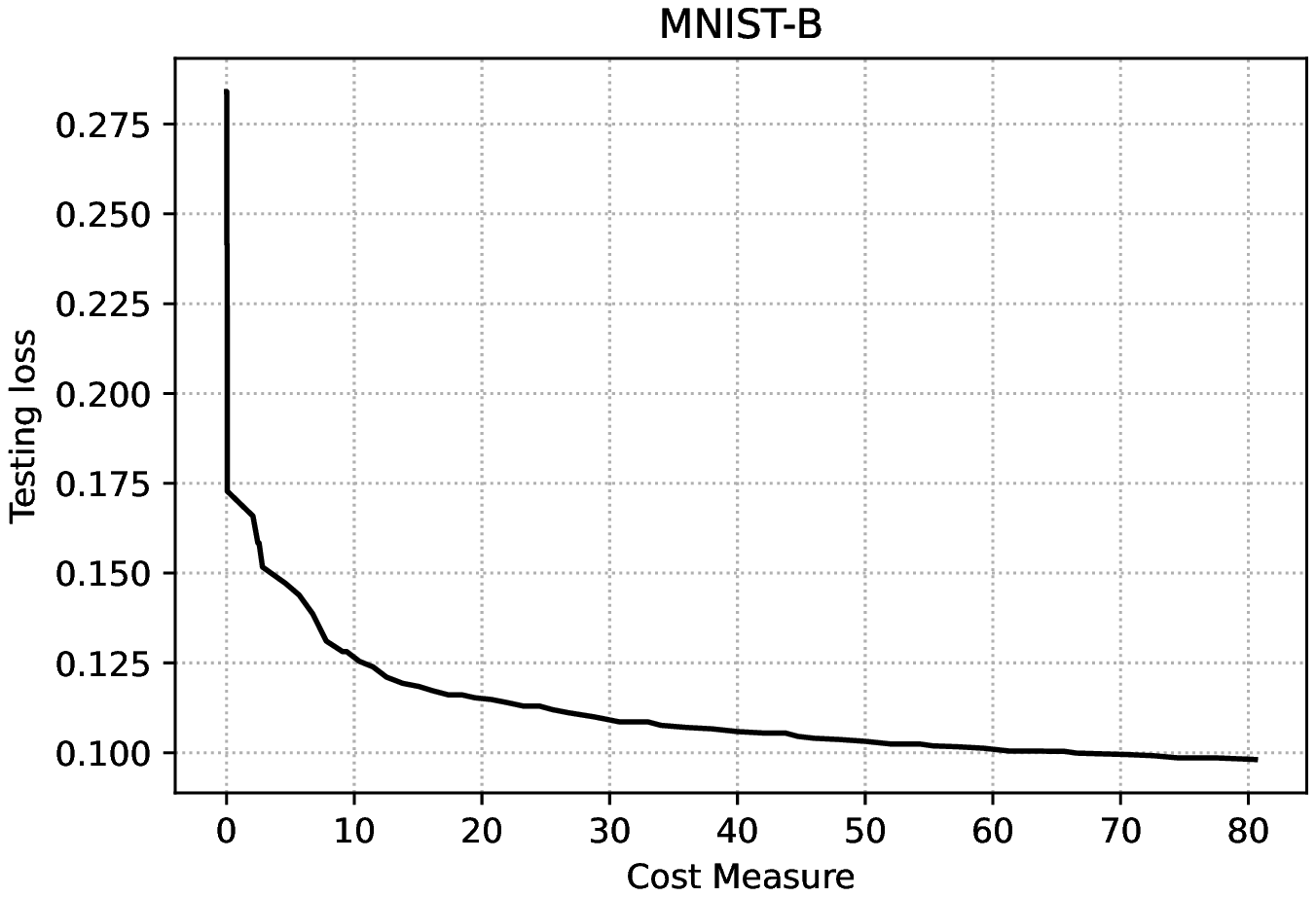}
\caption{MNIST-B  dataset. Training loss (left) and testing loss (right) by means of $IAR_1$ without net, against CM.}
\label{PerfAR1DA}
\end{figure}

 \begin{figure}[h]
\centering
\includegraphics[width=%
0.24\textwidth]{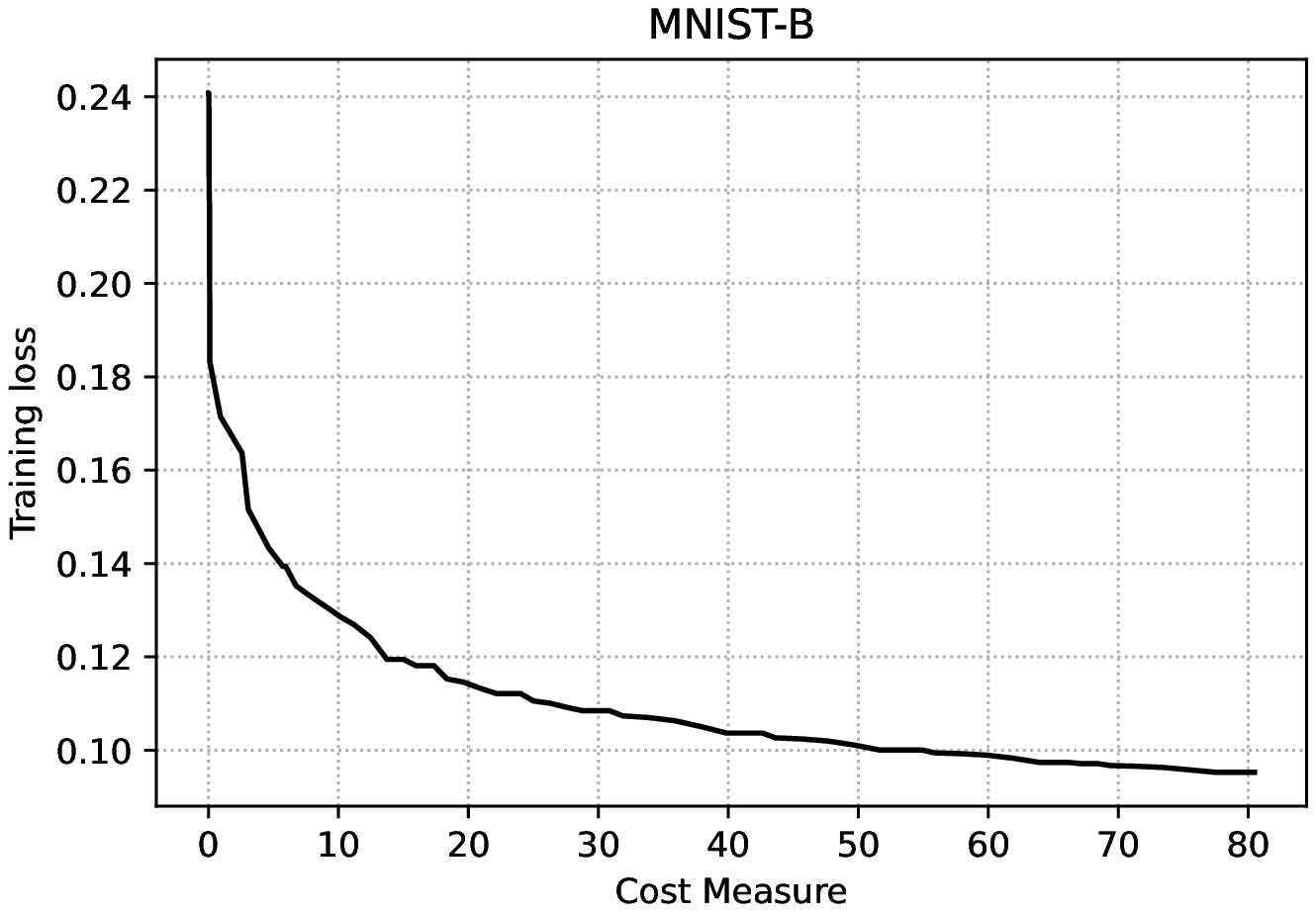}
\includegraphics[width=%
0.24\textwidth]{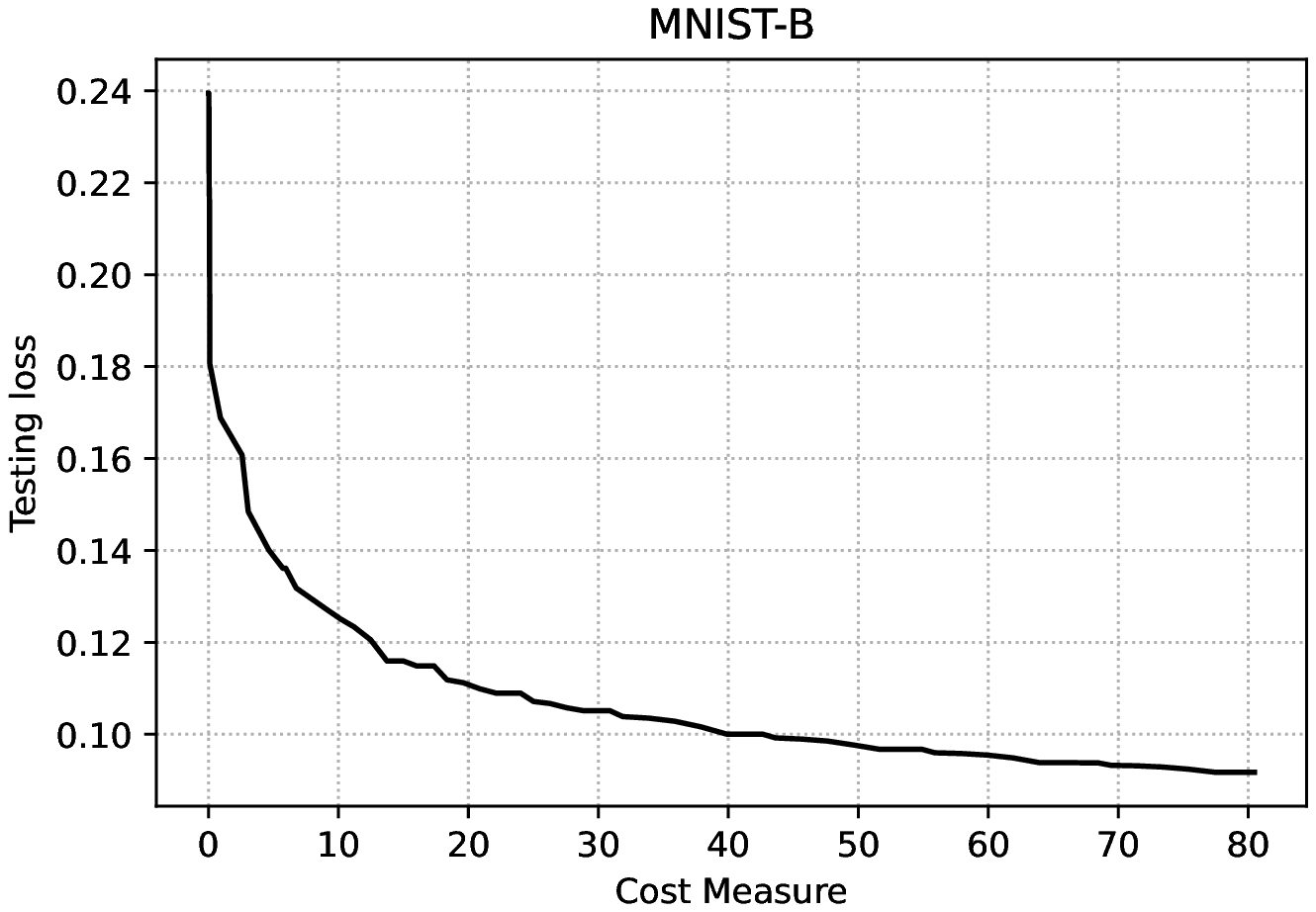}
\includegraphics[width=%
0.24\textwidth]{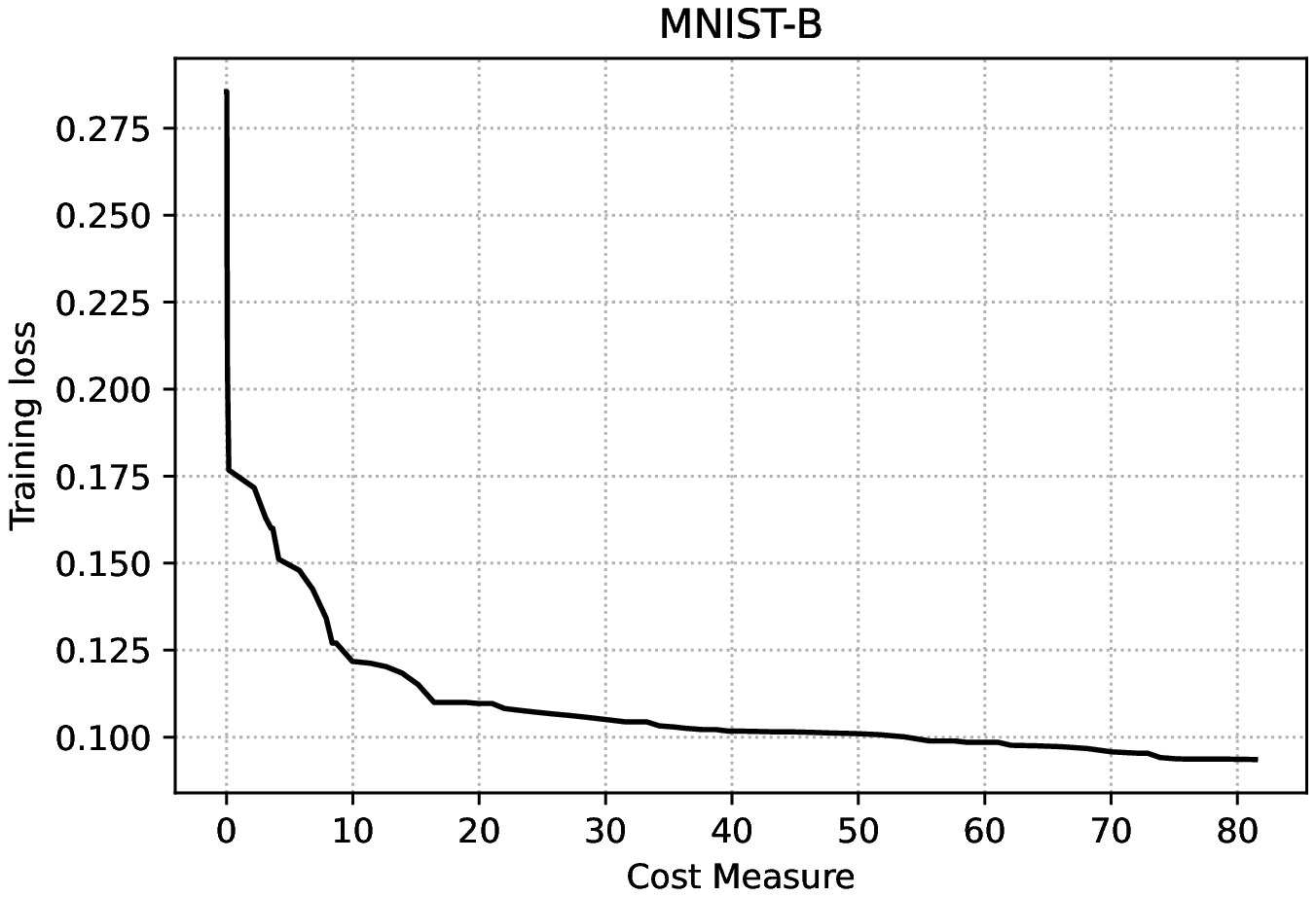}
\includegraphics[width=%
0.24\textwidth]{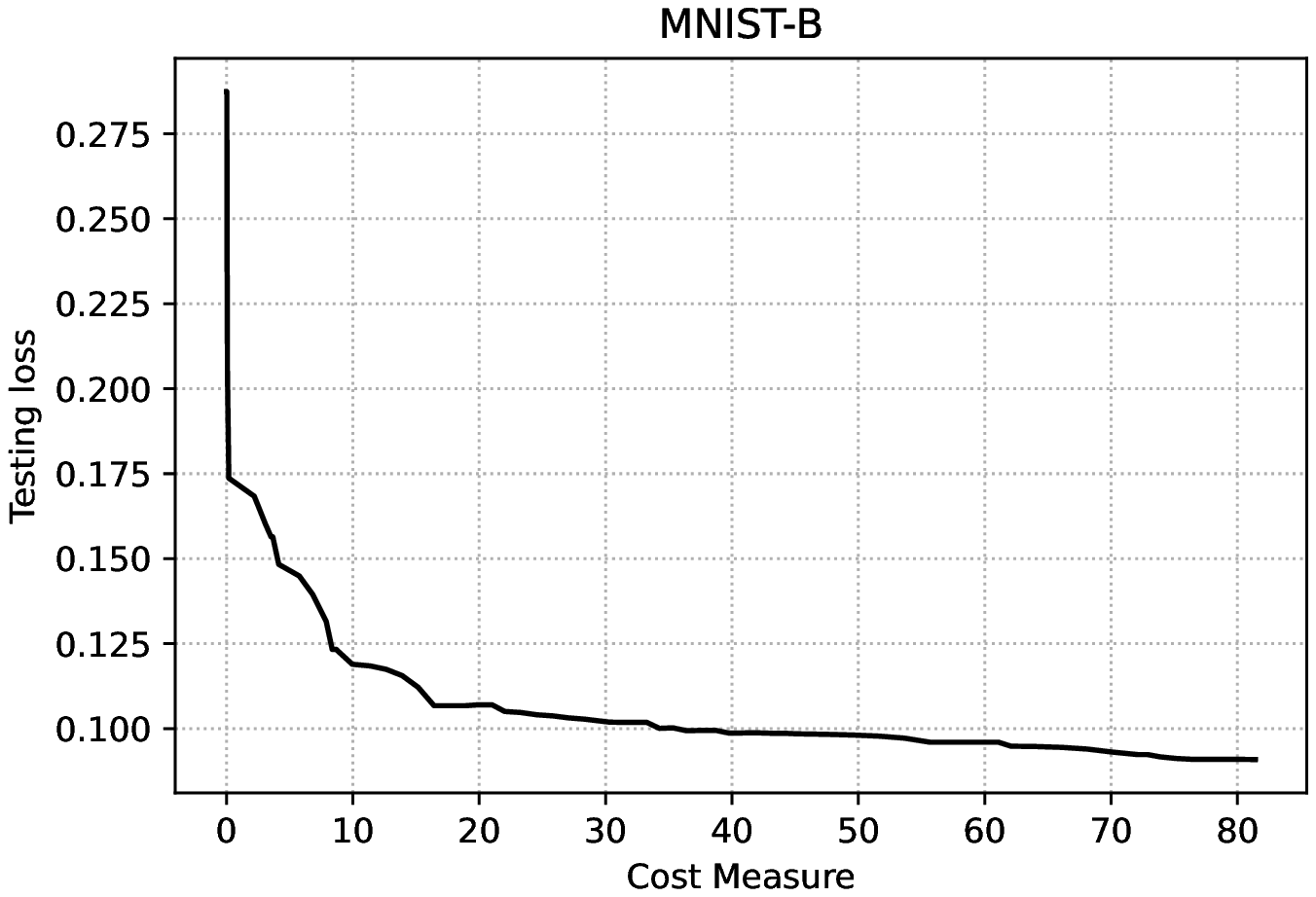}
\caption{MNIST-B dataset. Each row corresponds to the performance of $IAR_1$ with a different ANN architecture: shallow net $(15)$ (first row) and $2$-hidden layers net $(15,2)$ (second row). Training loss (left) and testing loss (right) against CM.}
\label{PerfANNAR1DA}
\end{figure}

Figures \ref{SSizeAR1DA}--\ref{SSizeANNAR1DA}, reporting the the adaptive sample sizes (in percentage) for computed function and gradient against CM within the $IAR_1$ algorithm, show that the sample size for estimating the objective function oscillates for a while, flattening on the full sample size eventually. The same is true for the corresponding gradients estimations, even if the growth toward the full sample is slower and the trend often remains below this threshold, especially for the MNIST-B dataset with ANN.

 \begin{figure}[h]
\centering
\includegraphics[width=%
0.24\textwidth]{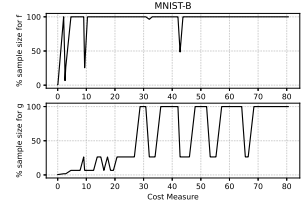}
\caption{MNIST-B dataset. Sample sizes used by $IAR_1$ with no net, against CM.}
\label{SSizeAR1DA}
\end{figure}

 \begin{figure}[h]
\centering
\includegraphics[width=%
0.24\textwidth]{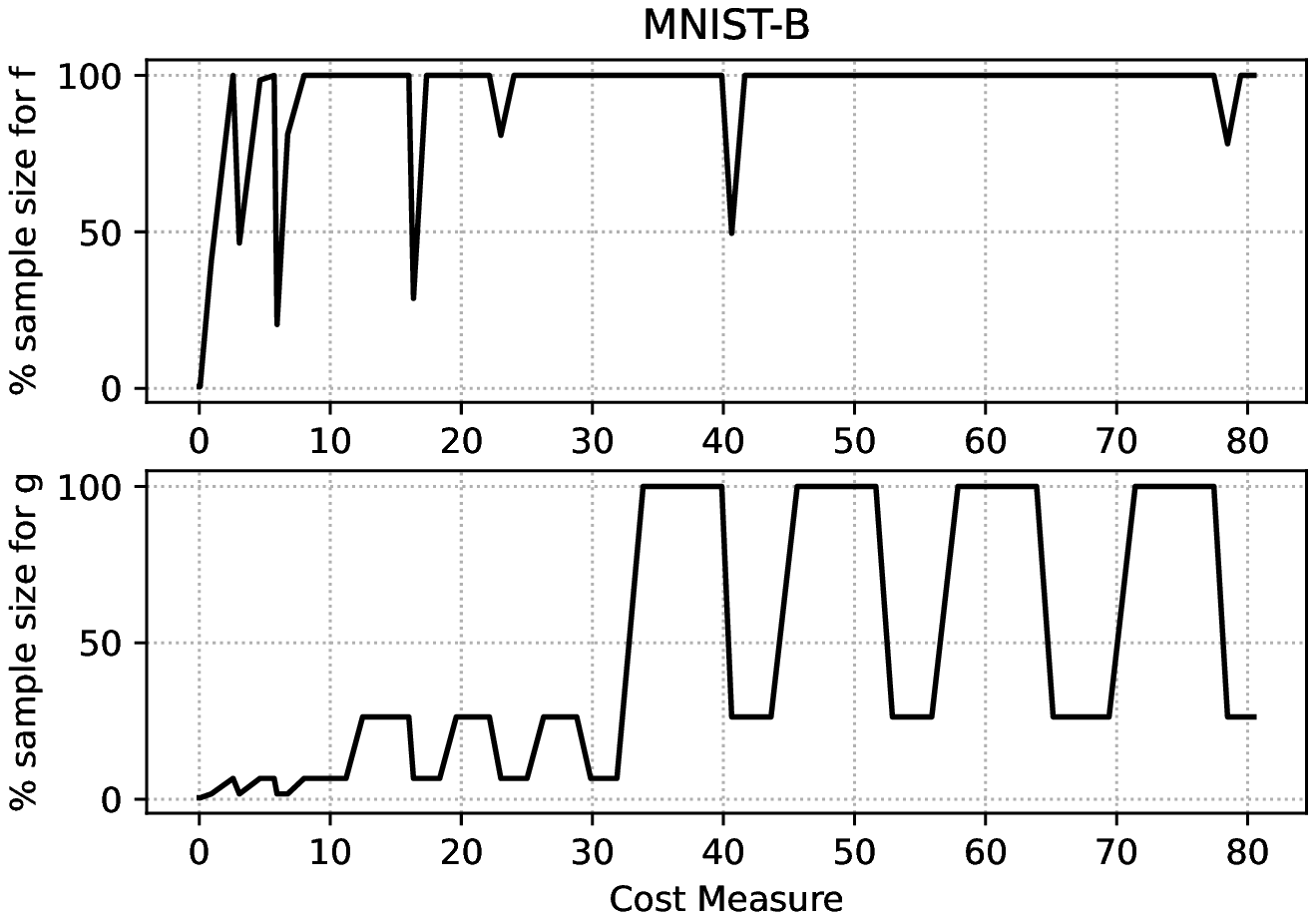}
\includegraphics[width=%
0.24\textwidth]{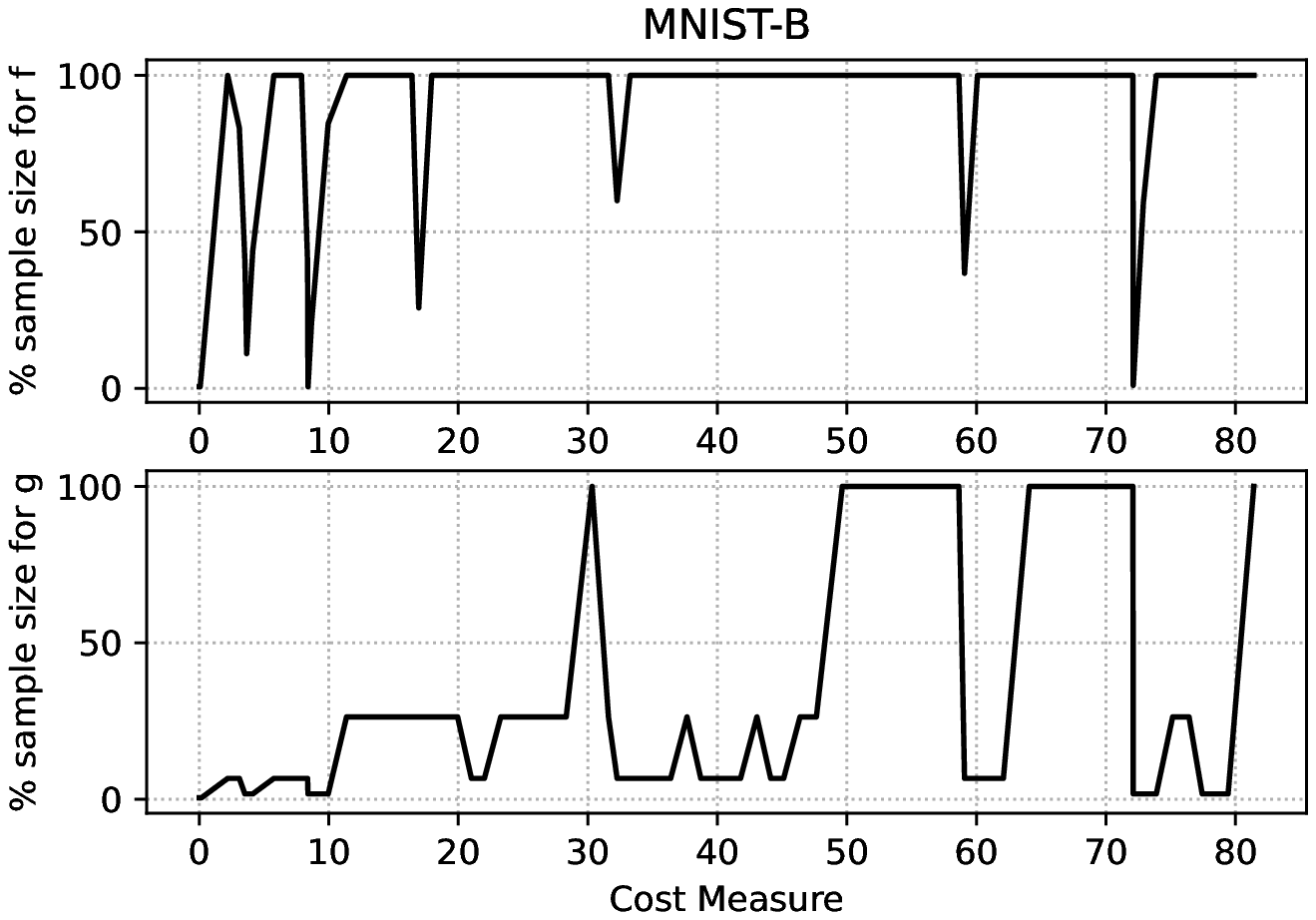}
\caption{MNIST-B dataset. Each column corresponds to the sample sizes used by $IAR_1$ with a different ANN architecture against CM: shallow net $(15)$ (left) and $2$-hidden layers net $(15,2)$ (right).}
\label{SSizeANNAR1DA}
\end{figure}
}

Moving to the GISETTE dataset,
{Figure \ref{PerfAR2DA} makes clear that the use of second order information enables the method to obtain lower values of the training and testing loss, with the same computational effort, and this yields  higher classification accuracy. We also note (see Figure \ref{SSizeAR2DA}) that $IAR_2$ employs  only occasionally the   $75\%$  of sample to approximate the Hessian matrix and in most of the iterations the Hessian approximation is obtained averaging on less than $20\%$ of samples. }

 \begin{figure}[h]
\centering
\includegraphics[width=%
0.44\textwidth]{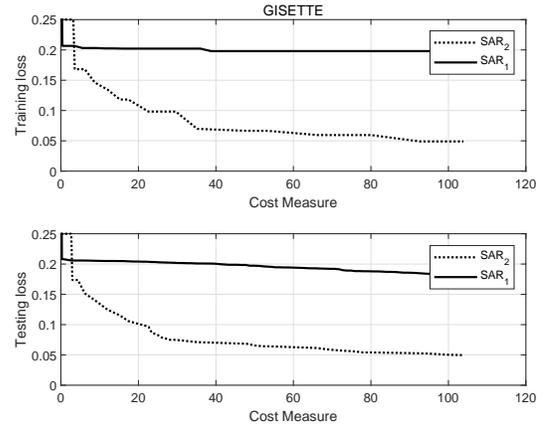}
\caption{GISETTE datasets. Training loss (top) and testing loss (bottom) by means of $IAR_1$ and $IAR_2$ without net, against CM.}
\label{PerfAR2DA}
\end{figure}

 \begin{figure}[h]
\centering
\includegraphics[width=%
0.44\textwidth]{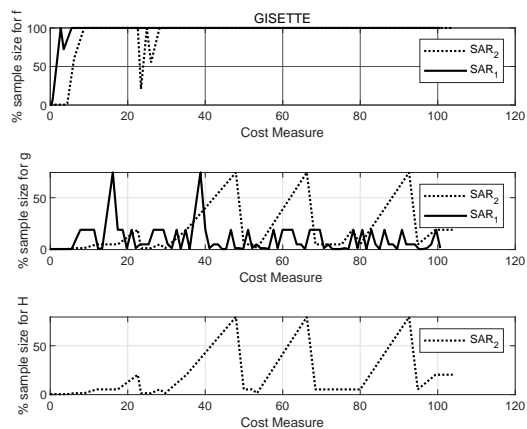}
\caption{GISETTE dataset. Sample sizes used by $IAR_1$ and $IAR_2$ against CM.}
\label{SSizeAR2DA}
\end{figure}

\section{Conclusions and perspectives}
The final expected bounds in \eqref{final-bound} are sharp in the order of the tolerances $\epsilon_1,\epsilon_2$. The effect of inaccurate evaluations is thus limited to scaling the
optimal complexity that we would otherwise derive from the deterministic analysis
(see, e.g., Theorem 5.5 in \cite{CartGoulToin20a} and Theorem 4.2 in \cite{IMA}), by a factor which depends on the probability $p_*$ of the model being accurate. Finally, first promising numerical results within the range of nonconvex finite-sum binary classification for imaging datasets confirm the expected behaviour of the method. From a theoretical perspective, the inclusion of inexact function evaluations subject to random noise is at the moment an open and challenging issue.
%
%
%
%

\section*{Acknowledgment}

INdAM-GNCS partially supported the first, second and third author under Progetti di Ricerca 2019 and 2020.

\end{document}